# ENUMERATION OF SOME CLOSED KNIGHT PATHS

**Stoyan Kapralov**
*Technical University of Gabrovo*

**Valentin Bakoev**
*V. Turnovo University*

**Kaloyan Kapralov**
*Skyscanner Bulgaria*

**Abstract**

*The aim of the paper is to enumerate all closed knight paths of length n over a square board of size n+1. The closed knight paths of length 4, 6 and 8 are classified up to equivalence. We determine that there are exactly 3 equivalence classes of closed knight paths of length 4, exactly 25 equivalence classes of closed knight paths of length 6 and exactly 478 equivalence classes of closed knight paths of length 8.*

## INTRODUCTION

The inspiration for this paper comes from the oldest lotto game in Bulgaria: "Toto 2" – 6 out of 49. The "Toto 2" game slip is a 7 x 7 grid pre-filled with the numbers from 1 to 49 starting with 1, 2, 3, … in the upper left and ending with … 47, 48, 49 in the lower right corner. Participants hoping to win the jackpot try to guess which 6 of these numbers will come up in a particular drawing.

In this paper, we consider those combinations which can be traversed with a closed path of the chess knight and we enumerate them up to equivalence.

Two combinations are equivalent if one can be obtained from the other by a sequence of one or more of the following transformations:
- translation - slide a combination over the board in one of the four directions (up, down, left, right);
- 90° rotation of the board;
- reflection of the board.

In fact the last two transformations are generators of a group of transformations which is isomorphic to $D_8$ - the dihedral group of order 8. The latter group is defined as a group of all symmetries of the square.

It is worth noting that the classic chess knight problem is the construction and classification of a *knight's tour*: a sequence of moves by a chess knight, which visits every square of the board exactly once. If from the last square visited, the knight can return in one knight's move to the starting square, the tour is a closed tour, otherwise it is an open tour. The original knight's tour problem is to find a closed tour on the classic 8 x 8 chessboard.

The first mathematical paper analyzing knight's tours was presented by the most productive mathematician of the eighteenth century, Leonhard Euler (1707–1783), to the Academy of Sciences at Berlin in 1759 (but not printed until 1766) [2].

Plenty of information about knight's tours is available on the internet, for example, see [3].

From a graph-theory point of view, finding a closed knight's tour is equivalent to finding a Hamiltonian cycle in the graph, corresponding to the chessboard. While it is well known that "Hamiltonian cycle" is an NP-complete problem, the special properties of the chessboard graph allow for the construction of a knight's tour in polynomial time [1, 6].

Although a single knight's tour can be constructed in polynomial time, finding the number of all knight's tours is a difficult problem, even in the computer age [7].

In 1996 two researchers write an algorithm to determine the total number of Knight's tours [4], yet their implementation turns out to be flawed within a year [5], when the total number of undirected tours is pegged at 13,267,364,410,532 with 1,658,420,855,433 equivalence classes under rotation and reflection of the board.

A knight's *tour* visits all squares of the board, while a knight's *path* may visit a subset of them, allowing for translation to also be considered as a transformation which produces equivalent paths.

In this paper we enumerate up to equivalence all closed paths of length $n = 4, 6$ and 8 on a square chessboard of size $n+1$.

**NEW RESULTS**

We determine that there are exactly 3 nonequivalent solutions for closed paths of length 4, exactly 25 nonequivalent solutions for closed paths of length 6, and exactly 478 nonequivalent solutions for closed paths of length 8 with a computer program whose algorithm is outlined below.

All nonequivalent solutions of length 4 are presented in Appendix A, and those of length 6 – in Appendix B.

The brief outline of our algorithm is as follows.

For each pair of cells we construct all connecting paths of length $k = n / 2$.

From every pair of paths we can construct a cycle of length $n$, if and only if, the paths do not intersect, that is, they have no common point except for the two endpoints.

It is sufficient to consider only cycles, which cannot be translated up or left, meaning that they have at least one point in the top row of the board and at least one point in the left-most column.

Each of these cycles can be framed in a rectangle. If the framing rectangles of two cycles are not congruent, then the cycles are not equivalent. Therefore, we only need to examine for equivalence between paths contained in each non-congruent type of rectangular frame.

**APPENDIX A**

```
1: 1 8 12 19
o . . . .
. . o . .
. o . . .
. . . o .

2: 2 9 11 18
. o . . .
```

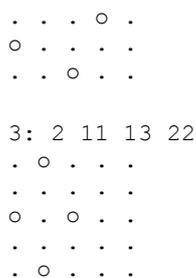

```
. . . o .
o . . . .
. . o . .

3: 2 11 13 22
. o . . .
. . . . .
o . o . .
. . . . .
. o . . .
```

**APPENDIX B**

```
 1:  1   2 10 11 15 16
 2:  1   3  5 10 16 18
 3:  1   3 10 12 16 25
 4:  1   3 10 16 18 23
 5:  1   5 10 11 16 20
 6:  1   5 10 16 18 31
 7:  1   5 10 16 20 25
 8:  1   6 10 11 16 19
 9:  1  10 11 16 19 24
10:  1  10 11 16 20 25
11:  1  10 16 19 25 34
12:  2   3  8 11 16 17
13:  2   4 15 17 19 24
14:  2  10 11 15 19 24
15:  2  10 11 15 20 25
16:  2  10 12 15 17 25
17:  2  10 15 17 19 32
18:  2  11 15 17 26 30
19:  2  11 15 20 24 33
20:  2  11 15 20 25 30
21:  2  11 15 24 26 39
22:  2  11 15 26 30 39
23:  2  15 17 22 24 37
24:  2  15 17 24 26 39
25:  2  15 17 30 32 45
```

**ACKNOWLEDGEMENTS**

The work of the first author was supported in part by Grant 1704C/2017 of the Technical University of Gabrovo, Bulgaria.

The second author is grateful for the partial support by the Research Fund of the University of Veliko Turnovo, Bulgaria, under Contract FSD-31-653-07/19.06.2017.